\documentclass[reqno,11pt]{amsart}
\usepackage{amsmath,amsthm,amsfonts,amssymb}
\usepackage[top=1in, bottom=1.25in, left=1.25in, right=1.25in]{geometry}

\usepackage{hyperref}
\usepackage[all]{xy}
\usepackage{graphicx}
\usepackage{amscd}
\usepackage{color}
\usepackage{enumerate,enumitem}
\usepackage{fourier}
\usepackage{tikz}
\usepackage{cases}

\hypersetup{
    colorlinks=true,
    linkcolor=blue,
    citecolor=green,
    urlcolor=cyan
}

\numberwithin{equation}{section}

\theoremstyle{plain} 
\newtheorem{theorem}{Theorem}[section]

\newtheorem{prop}[theorem]{Proposition}

\theoremstyle{definition}
\newtheorem{definition}[theorem]{Definition}

\newtheorem{thmx}{Theorem}

\theoremstyle{remark}

\newcommand{\PP}{\mathbb{P}}

\newcommand{\Z}{\mathbb{Z}}

\newcommand{\C}{\mathbb{C}}

\newcommand{\CO}{{\mathcal{CO}}}

\newcommand{\CB}{{\mathcal{B}}}

\newcommand{\ks}{\mathfrak{ks}}

\newcommand{\AI}{{A_\infty}}

\newcommand{\Jac}{{\rm Jac}}

\newcommand{\pt}{{\rm{pt}}}
\newcommand{\HG}{{\widehat{G}}}

\newcommand{\bL}{{\mathbb{L}}}

\newcommand{\barW}{{\overline{W}}}
\newcommand{\tildeW}{{\widetilde{W}}}

\newcommand{\Hess}{{\rm{Hess}}}

\newcommand{\WT}[1]{\widetilde{#1}}

\allowdisplaybreaks

\begin{document}
\title[Kodaira-Spencer maps and Frobenius algebras]
{Kodaira-Spencer maps for elliptic orbispheres \\ as isomorphisms of Frobenius algebras}

\author{Sangwook Lee}
\address{Sangwook Lee: Department of Mathematics and Integrative Institute of Basic Science \\Soongsil University \\
369 Sangdo-ro, Dongjak-gu, Seoul, Korea}
\email{sangwook@ssu.ac.kr}
\begin{abstract}
Given a mirror pair of a symplectic manifold $X$ and a Landau-Ginzburg potential $W$, we are interested in the problem whether the quantum cohomology of $X$ and the Jacobian algebra of $W$ are isomorphic. Since those can be equipped with Frobenius algebra structures, we might ask whether they are isomorphic as Frobenius algebras. We show that the Kodaira-Spencer map gives a Frobenius algebra isomorphism for elliptic orbispheres, under the Floer theoretic modification of the residue pairing.
\end{abstract}
\maketitle

\section{Introduction}
Closed string mirror symmetry predicts that quantum cohomology of a symplectic manifold and Jacobian ring of the mirror superpotential are isomorphic. There have been results on this problem firstly for toric manifolds: see \cite{batyrev93,givental,iritani09} etc.

Fukaya-Oh-Ohta-Ono also gave a construction of the ring isomorphism in \cite{FOOO-T,FOOO-T2,FOOO10b}, based on the study of closed-open map which is defined more geometrically. Though the codomain of the original closed-open map is Hochschild cohomology which is hard to grasp, Fukaya-Oh-Ohta-Ono proved that the length-0 part of the closed-open map (which was called {\em Kodaira-Spencer map}), whose codomain is now Jacobian ring, provides an isomorphism in compact toric case. Their strategy was also employed by Amorim-Cho-Hong-Lau in \cite{ACHL} for proving mirror symmetry for orbifold projective lines.

Since quantum cohomology ring and Jacobian ring are both Frobenius algebras with respect to the Poincar\'e duality and the residue pairing respectively, we might ask whether their pairings are also related in a suitable way. The question was dealt with in \cite{CLS}, and it was conjectured that if we want the Kodaira-Spencer map to be a Frobenius algebra isomorphism, we need to modify the residue pairing by a constant which is the ratio of "Floer volume form" and the usual volume form on a Lagrangian submanifold. The rescaling constant appeared naturally when we consider Cardy condition (see \cite{CLS} for more detail). 

In this paper we focus on elliptic orbispheres. Several mirror isomorphisms of Frobenius manifolds for orbifold projective lines have been established. In \cite{SaTa} Satake-Takahashi proved that there is an isomorphism of Frobenius manifolds from Gromov-Witten theory of $\PP^1_{3,3,3}$ and the universal unfolding of the mirror potential, and also for $\PP^1_{2,2,2,2}$ case (where the counterpart is given by the invariant theory of an elliptic Weyl group). For spherical cases, namely for $\PP^1_{a,b,c}$ with $\frac{1}{a}+\frac{1}{b}+\frac{1}{c}>1$ we refer readers to \cite{Ros}, \cite{milanov-tseng08}, \cite{IST} etc. Towards the ultimate level of closed string mirror symmetry such as above results, we add more elliptic orbisphere examples, namely $\PP^1_{2,4,4}$ and $\PP^1_{2,3,6}$ which are quotients of the elliptic curve by $\Z_4$ and $\Z_6$ respectively. $\PP^1_{3,3,3}$ will be also revisited. Though we work only at the level of Frobenius algebras, we hope that we can find a relationship with previous works, such as the relation between the Floer theoretic rescaling constant and the choice of primitive form.

We summarize the idea. Let $W$ be the mirror superpotential to $X$, so the Kodaira-Spencer map is given by $\ks: QH^*(X) \to \Jac(W)$. The most natural pairings for Frobenius algebra structure are Poincar\'e pairing and the residue pairing respectively. We point out that the residue pairing $\langle\;,\;\rangle_{res}$ on $\Jac(W)$ will be chosen as
\begin{equation}\label{eq:residue}
\langle f,g\rangle_{res}=(-1)^{\frac{n(n-1)}{2}}{\rm{Res}}\begin{bmatrix}
        fg\cdot dx_1\wedge\cdots\wedge dx_n \\ \partial_{x_1}W,\cdots,\partial_{x_n}W
    \end{bmatrix}\end{equation}
whose formula appears in \cite{PV} as induced from the Mukai pairing on $HH_*(MF(W))$, so that it is "the most canonical" in some sense. It is also remarkable that the sign $(-1)^{\frac{n(n-1)}{2}}$ also appears in the Cardy condition in \cite[Theorem 3.4.1]{FOOO10b}.

Our main theorem is as follows.
\begin{thmx}\label{thm:main}
    Let $X$ be an elliptic orbisphere and $W$ be its mirror superpotential. Let $\bL \subset X$ be the Seidel Lagrangian with odd degree immersed generators $X,Y,Z$, and $c_\bL$ be the constant defined by 
\begin{equation}\label{eq:volconst} c_\bL\cdot p = {m_2(X,m_2(Y,Z))}\end{equation}
 where $p =m_2(X,\bar{X})$.
 Then the Kodaira-Spencer map \[\ks: \big(QH^*(X),\langle \cdot,\cdot\rangle_{PD}\big) 
    \to \big(\Jac(W),\langle c_\bL \cdot,c_\bL\cdot \rangle_{res}\big)\] is an isomorphism of Frobenius algebras.
\end{thmx}

Observe that we modified the residue pairing by suitable constant $c_\bL$ as discussed above. To compare pairings we need to compute the residue of $\ks(\pt_X)$, where $\pt_X$ is the Poincar\'e dual of point class of the symplectic manifold $X$. Though $\ks(\pt_X)$ is explicitly computed in \cite{ACHL} for general orbifold sphere $X$, it is hard to conclude that its (rescaled) residue is indeed $1=\int_X \pt_X$ as expected. The difficulty arises in the comparison of two different arithmetics of formal power series. 
We bypass this difficulty by computing $\ks(\pt_X)$ in another way. Let $G$ act on the elliptic curve $E$ so that $X=E/G$. We will use the result of \cite{CLe} that the Kodaira-Spencer map on $X$ can be lifted to the {\em orbifold} Kodaira-Spencer map $\ks_{orb}:H^*(E)\to \Jac(W,\HG)$, where $\Jac(W,\HG)$
is the {\em orbifold} Jacobian algebra ($\HG$ is the character group of $G$, so it is isomorphic to $G$ if $G$ is abelian). 
By the relation
\[ \ks(\pt_X)=\frac{1}{|G|}\ks_{orb}(\pt_E),\]
we can recover $\ks(\pt_X)$ by $\ks_{orb}(\alpha\cup\beta)=\ks_{orb}(\alpha)\bullet\ks_{orb}(\beta)$ for some $\alpha,\beta\in H^1(E)$. By a classical result on the residue over an isolated singularity (which will be recalled in Section \ref{sec:mainthm}), the following is a rephrasing of Theorem \ref{thm:main}.

\begin{thmx}\label{thm:main'}
Let $X$ be an elliptic orbisphere $E/G$ and $W$ be its mirror superpotential. 
 Then
 \begin{equation}\label{eq:main}
 c_\bL^2\cdot\ks_{orb}(\pt_E)=-\frac{\det \; \Hess(W)}{|G|\cdot \mu}.\end{equation}
\end{thmx}

The organization of the paper is as follows. In Section \ref{sec:KS} we first recall the construction of orbifold Jacobian algebras from Floer theory. Then we briefly review (orbifold) Kodaira-Spencer maps which appear in Fukaya-Oh-Ohta-Ono's works and also in \cite{CLe}. In Section \ref{sec:orbiproduct} we review the product structure on orbifold Jacobian algebras following \cite{LeeTwJac}, and explicitly compute the product for the mirror Landau-Ginzburg orbifold to an elliptic curve. Finally, we prove our main result in Section \ref{sec:mainthm}, with a remark on nontrivial identities of arithmetics of formal power series.

\begin{center}
    {\bf Acknowledgement}
\end{center}

This work was supported by Basic Science Research Program through the National Research Foundation of Korea (NRF) funded by the Ministry of Education (2021R1A6A1A10044154).

\section{Kodaira-Spencer map to an orbifold Jacobian algebra}\label{sec:KS}
\subsection{Preliminaries}
We briefly recall the deformation theory of Lagrangian submanifolds originally in \cite{FOOO}, following every notation in Section 4.1 of \cite{CLe}. Let $X$ be a symplectic manifold and $\bL$ be its (possibly immersed) Lagrangian submanifold. Let $CF(\bL,\bL;\Lambda)$ be the Fukaya $\AI$-algebra whose underlying space is the sum of de Rham algebra and the module generated by immersed generators. Suppose that $\bL$ is weakly unobstructed, i.e. the $\AI$-structure on $CF^*(\bL,\bL;\Lambda)$ can be deformed by weak Maurer-Cartan elements. Assume further that $X_1,\cdots,X_n\in CF^1(\bL,\bL;\Lambda)$ are weak Maurer-Cartan elements such that $c_1 X_1+\cdots +c_n X_n$ is also a weak Maurer-Cartan element for any $c_i\in \Lambda_+$. Let $x_i$ be the dual variable of $X_i$ and $b:=x_1 X_1+\cdots+x_n X_n$. The weak Maurer-Cartan equation gives rise to the potential $W_\bL\in R=\Lambda[x_1,\cdots,x_n]$. We consider the following $\AI$-algebra
\[ \CB(\bL):= (CF^*(\bL,\bL;\Lambda)\otimes_{\Lambda_0} R, \{m_k^b\}). \]
We also denote by $\CB(\bL)_{alg}$ the associative algebra with the same underlying space as $\CB(\bL)$ while equipped with the product $v\cdot w:=(-1)^{|v|}m_2^b(v,w)$.
The following results evidently reflect the importance of the above definition.
\begin{prop}[\cite{CLe}]
    If $\bL$ is the Seidel Lagrangian in an orbisphere $\PP^1_{a,b,c}$ or the Lagrangian torus at the critical point of $W_\bL$, then there is an algebra isomorphism 
    \begin{equation}\label{eq:FloertoJac}
    \Psi: H^*(\CB(\bL)_{alg}) \stackrel{\simeq}{\longrightarrow} \Jac(W_\bL)=R/\partial W_\bL.\end{equation}
\end{prop}
Let us now recall the orbifold Jacobian algebra of an isolated singularity equipped with a group action. 
\begin{definition}
 Let $H$ be a finite abelian group, which acts on $R=\Bbbk[x_1,\cdots,x_n]$ and leaves $W$ invariant. We call the pair $(W,H)$ a {\em Landau-Ginzburg orbifold}.
\end{definition}
Throughout the paper, we only consider diagonal $H$-action, i.e. $h\cdot x_i=h_i x_i$ for some $h_i \in \Bbbk^*$.
\begin{definition}
Let $(W,H)$ be a Landau-Ginzburg orbifold. Then the {\em twisted Jacobian algebra} of $(W,H)$ is defined as
\[ \Jac'(W,H):= \bigoplus_{h\in H} \Jac(W^h)\cdot \xi_h\]
where $W^h$ is the image of $W$ via projection $\pi: R \to R/(x_i: hx_i \neq x_i)$. The formal generator $\xi_h$ has degree $|I_h|\in \Z_2$ where $I_h=\{i\in \{1,\cdots,n\} \mid hx_i\neq x_i\}$. The $H$-action on generators is defined by
\[ h'\cdot \xi_h = \frac{\prod_{i\in I_h} (h'^{-1}x_i)}{\prod_{i\in I_h}x_i} \xi_h,\]
and the $H$-invariant subalgebra \[\Jac(W,H):=\Jac'(W,H)^H\] is called the {\em orbifold Jacobian algebra} of $(W,H)$.
\end{definition}
\begin{theorem}
\begin{enumerate}
    \item $\Jac'(W,H)$ is an $H$-graded algebra, namely \[\xi_h \bullet \xi_{h'} \in \Jac(W^h)\cdot \xi_{hh'}.\]
    \item The product $\bullet$ satisfies braided-commutativity, namely
    \[ \xi_h \bullet \xi_{h'} = (-1)^{|\xi_h|\cdot |\xi_{h'}|}h\cdot \xi_{h'} \bullet \xi_h.\]
    In particular, if $\xi_h$ and $\xi_{h'}$ are both $H$-invariant, then 
    \[ \xi_h \bullet \xi_{h'} = (-1)^{|\xi_h|\cdot |\xi_{h'}|} \xi_{h'} \bullet \xi_h,\]
    which implies that $\Jac(W,H)$ is a supercommutative algebra.
\end{enumerate}

\end{theorem}
We postpone the definition of the product on $\Jac(W,H)$ until Section \ref{sec:orbiproduct}.
Orbifold Jacobian algebras appear naturally in Floer theory as follows. Let $X$ be a symplectic manifold and a finite abelian group $G$ act on $X$. Suppose that $\bL \subset X/G$ is a weakly unobstructed Lagrangian submanifold with mirror superpotential $W$, with an embedded Lagrangian lift $\bL_0\subset X$. Then we can construct a new $\AI$-algebra structure on $\CB(\bL)\otimes \Lambda[\HG]$, whose $\HG$-invariant subalgebra $(\CB(\bL)\otimes \Lambda[\HG])_{alg}^\HG$ is isomorphic to $\Jac(W,\HG)$. 

\begin{prop}[\cite{CLe}]\label{prop:uptodown}
    Let $\widetilde{\bL}=\bigoplus_{g\in G}g\cdot \bL_0\subset X$ be a Lagrangian submanifold with lifted weak bounding cochain $\WT{b}$ from $b$. Then
    \[ \Phi: \big(CF(\widetilde{\bL},\widetilde{\bL})\otimes R, m_k^{\WT{b}}\big) \to \big((\CB(\bL)\otimes \Lambda[\HG])^\HG,m_k^{b\otimes 1}\big),\; v_g \mapsto \frac{1}{|G|}\sum_{\chi\in \HG}\chi(g)v\otimes \chi\]
    is an $\AI$-isomorphism. Here $v_g$ is an element in $CF(g\cdot \bL_0, h\cdot \bL_0)$ and $v$ is the projection of $v_g$.
\end{prop}
Despite being vacuous, $\Phi$ can be also defined in the previous nonequivariant setting (namely, endowed with trivial group action) just by the identity.

\begin{prop}[\cite{CLe}]\label{prop:FloertoorbJac}
    Let $(\CB(\bL)\otimes \Lambda[\HG])^\HG_{alg}$ be an associative algebra with 
    \[(v\otimes \chi)\cdot(w\otimes \chi'):=(-1)^{|v|}m_2^{b\otimes 1}(v\otimes \chi,w\otimes \chi').\] Then there is an algebra isomorphism
    \begin{equation}\label{eq:FloertoorbJac}
    \Psi_\HG:H^*((\CB(\bL)\otimes \Lambda[\HG])^\HG_{alg}) \stackrel{\simeq}{\longrightarrow} \Jac(W,\HG).
    \end{equation}
\end{prop}
Again, in the nonequivariant setting, $\Psi_\HG$ is nothing but $\Psi$ above.

Now we recall the construction of Kodaira-Spencer maps. For general weakly unobstructed Lagrangian $L\subset X$ with bounding cochain $b$, we use the similar configuration of holomorphic discs as Fukaya-Oh-Ohta-Ono's original definition. 
\begin{definition}
For a moduli space of $J$-holomorphic discs with $k+1$ boundary marked points and one interior marked points $\mathcal{M}_{k+1,1}(\beta)$ for $\beta \in H_2(X,L)$ and a cycle $A \subset X$,
let $\mathcal{M}_{k+1,1}(\beta,A) = \mathcal{M}_{k+1,1}(\beta) \times_X A$ and consider their evaluation maps $ev_i^k:\mathcal{M}_{k+1,1}(\beta,A) \to L$ at the $i$th marked point.
Consider the length-0 closed open map
\begin{equation}\label{eq:defks}
\CO^0: QH^*(X) \to H^*(\CB(L)_{alg}), \; \CO^0 (PD[A]) := \sum_{\beta\in H_2(X,L)} \sum_{k=0}^\infty (ev_0^k)_! ((ev_1^k)^*b \wedge \cdots \wedge (ev_k^k)^*b).
\end{equation}
\end{definition}

\begin{definition}
Let $X$ be a symplectic manifold with finite group $G$ acting effectively on it($G$ may be trivial), and $L$ be a $G$-equivariant weakly unobstructed Lagrangian. Then the {\em orbifold Kodaira-Spencer map} is
$\ks_{orb}:=\Psi_\HG \circ \Phi \circ \CO^0.$
\end{definition}
Note that if $G$ is trivial, then $\ks_{orb}$ is just the ordinary Kodaira-Spencer map $\ks$. It was shown in \cite{CLe} that $\ks_{orb}$ is well-defined, and furthermore it is a ring homomorphism.

\subsection{Computation of orbifold Kodaira-Spencer maps from elliptic curves}\label{sec:KSelliptic}
We choose three different Lagrangian submanifolds on elliptic curves and compute the orbifold Kodaira-Spencer maps with respect to them. 
\subsubsection{$\Z_3$}
\begin{figure}
    \centering
    \includegraphics[width=0.85\linewidth]{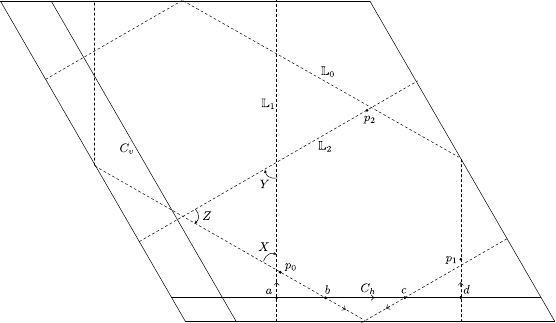}
    \caption{$\Z_3$-equivariant Lagrangian}
    \label{fig:333}
\end{figure}

Let $E= \C / (\Z + e^{2\pi i /3}\Z)$ be an elliptic curve and let $\Z_3=\{1,\rho=e^{2\pi i/3},\rho^2=e^{4\pi i/3}\}$ act on $E$ by multiplication, so $\PP^1_{3,3,3}=E/\Z_3$. Let $\bL_0 \subset E$ be an embedded circle in $E$, and $\bL_1:=\rho \bL_0$, $\bL_2:=\rho^2 \bL_0$ as in Figure \ref{fig:333}. Then $\WT{\bL}:=\bL_0\oplus\bL_1\oplus\bL_2$ is a weakly unobstructed Lagrangian on $E$ with potential $W_{333}$.

Let $C_h$ be a homology cycle representing a class $(1,0)$. 
Since $\dim\mathcal{M}_{k+1,1}(\beta, C_h)=k+\mu(\beta)$, if $\beta\in H_2(E,\WT{\bL})$ is nontrivial then there is no summand in \eqref{eq:defks} for $\beta$. Therefore to compute the Kodaira-Spencer map we only consider $\mathcal{M}_{1,1}(0,C_h)$. There is a natural orientation on $\mathcal{M}_{1,1}(0)$ because it is diffeomorphic to $\WT{\bL}$, and the fiber product $\mathcal{M}_{1,1}(0,C_h)$ is nothing but the intersection of $\WT{\bL}$ and $C_h$. We conclude that $\ks(PD[C_h])$ is given by (Poincar\'e dual of) the oriented intersection $\WT{\bL}\cap C_h$. In Figure \ref{fig:333} we depicted intersection points between $\WT{\bL}$ and $C_h$ whose Poincar\'e dual 1-forms are $a,b,c,d$ respectively.
Taking orientations into account, we have 
\[ \CO^0(PD[C_h])= -a+b+c-d.\]

To read an orbifold Jacobian algebra element from it, we need to recall the construction of an isomorphism $H^*\big((\CB(\bL)\otimes \Lambda[{\widehat{\Z_3}}])^{\widehat{\Z_3}}\big) \cong \Jac(W_{333},{\widehat{\Z_3}})$ for the Seidel Lagrangian $\bL$. 
The module $CF(\bL,\bL)$ is generated by $\big\{1,X,Y,Z,X \wedge Y,Y\wedge Z,Z\wedge X,X\wedge (Y\wedge Z)\big\}$ where $X,Y$ and $Z$ are odd degree immersed generators and $\wedge$ is the binary $\AI$-product $m_2$ without weak bounding cochain insertions. 
Observe that 
\[ \frac{X\wedge Y}{\bar{Z}}=\frac{Y\wedge Z}{\bar{X}}=\frac{Z\wedge X}{\bar{Y}}=\frac{X\wedge (Y\wedge Z)}{p}=c_\bL\]
for some constant $c_\bL$, as in \eqref{eq:volconst} (recall that $p=m_2(X,\bar{X})$).
\begin{theorem}[\cite{CLe}]
The isomorphism 
\[\Psi_{\widehat{\Z_3}}:H^*(\CB(\bL)\otimes \Lambda[\widehat{\Z_3}])^{\widehat{\Z_3}} \cong \Jac(W_{333},\widehat{\Z_3})\] in \eqref{eq:FloertoorbJac} is given by
\begin{align}
\begin{split}
    H^*(\CB(\bL)\otimes 1)^{\widehat{\Z_3}} &\to \Jac(W_{333})^{\widehat{\Z_3}}, \qquad\quad\;\; f\cdot 1 \mapsto f, \\
    H^*(\CB(\bL)\otimes \chi)^{\widehat{\Z_3}} &\to \Jac(W_{333}^\chi)^{\widehat{\Z_3}}\cdot \xi_\chi, \quad \big(X\wedge(Y\wedge Z)+ ({\rm lower})\big)\otimes \chi \mapsto \xi_\chi. \label{eq:BtoJac}
    \end{split}
\end{align} 
\end{theorem}
For $p=m_2(X,\bar{X}) \in CF(\bL,\bL)$, let $p_i\in CF(\bL_i,\bL_i)$ for $i=0,1,2$ such that $p_i$ projects to $p$. Now let us investigate the Poincar\'e dual of each intersection of $\WT{\bL}$ and $C_h$, say $a$ for example. Seeing Figure \ref{fig:333} again, $p_1$ and $a$ are cohomologous in the de Rham complex $\Omega^*(\bL_1)$, but they are not cohomologous in $(CF(\WT{\bL},\WT{\bL}),m_1^{\WT{b}})$ (in fact, they are not even cocycles). If we consider a de Rham 0-form $I$ whose de Rham coboundary is $p_1-a$, then $m_1^{\WT{b}}(I)=p_1-a+{\rm (lower)}$, where (lower) means a linear sum of odd degree immersed generators. In the same vein, we consider de Rham 0-forms $J$, $K$ and $L$ whose coboundaries are $p_0-b$, $p_2-c$ and $p_1-d$ respectively. Then
\[m_1^{\WT{b}}(I-J-K+L)=(-a+b+c-d)+(p_1-p_0-p_2+p_1)+{\rm (lower)},\]
i.e. $-a+b+c-d$ is cohomologous to $-2p_1+p_0+p_2+{\rm (lower)}$ in $CF(\WT{\bL},\WT{\bL})$ with respect to $m_1^{\WT{b}}$. Therefore, letting $\widehat{\Z_3}=\{1,\chi,\chi^2\},$ the image via (orbifold) Kodaira-Spencer map is
\begin{align*} \Phi(\CO^0(PD[C_h]))&=\frac{1}{3} \big( (-2\chi(\rho)+1+\chi(\rho^2))p\otimes \chi + (-2\chi^2(\rho)+1+\chi^2(\rho^2))p\otimes \chi^2\big)+{\rm (lower)}\\
&\stackrel{\eqref{eq:BtoJac}}{\mapsto}\displaystyle\frac{(-2\chi(\rho)+1+\chi(\rho^2))\xi_\chi+(-2\chi^2(\rho)+1+\chi^2(\rho^2))\xi_{\chi^2}}{3c_\bL}\in \Jac(W_{333},{\widehat{\Z_3}}).
\end{align*}
Observe that there is no output on $1$-sector due to degree reason. For the cycle $C_v$ of class $(0,1)$,
\begin{align*} \Phi(\CO^0(PD[C_v]))&=\frac{1}{3}\big( (-2\chi(\rho^2)+1+\chi(\rho))p\otimes \chi + (-2\chi^2(\rho^2)+1+\chi^2(\rho))p\otimes \chi^2\big)+{\rm (lower)}\\
&\stackrel{\eqref{eq:BtoJac}}{\mapsto}\displaystyle\frac{(-2\chi(\rho^2)+1+\chi(\rho))\xi_\chi+(-2\chi^2(\rho^2)+1+\chi^2(\rho))\xi_{\chi^2}}{3c_\bL}\in \Jac(W_{333},{\widehat{\Z_3}}).
\end{align*}
Letting $\chi(\rho)=\rho$, we summarize
\begin{equation}\label{eq:ks333} 
    \ks_{orb}(PD[C_h])= -\frac{\rho}{c_\bL}\xi_\chi-\frac{\rho^2}{c_\bL}\xi_{\chi^2},\quad
    \ks_{orb}(PD[C_v])= -\frac{\rho^2}{c_\bL}\xi_\chi-\frac{\rho}{c_\bL}\xi_{\chi^2}.
\end{equation}
We hope readers notice that $c_\bL$ is involved in the computation.

\subsubsection{$\Z_4$}
\begin{figure}
    \centering
    \includegraphics[width=0.4\linewidth]{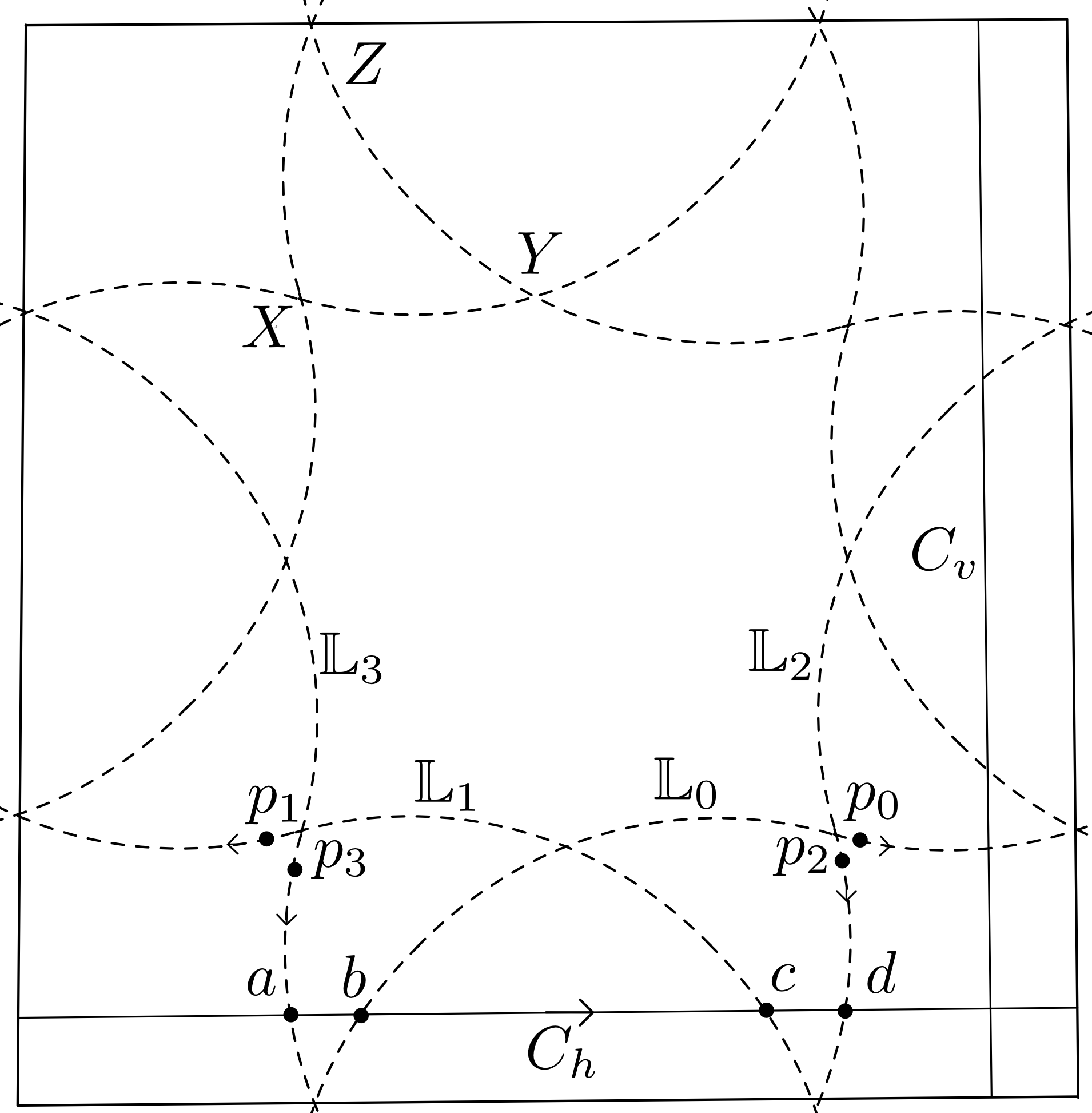}
    \caption{$\Z_4$-equivariant Lagrangian}
    \label{fig:244}
\end{figure}

Let $E=\C/(\Z + i\Z)$ be an elliptic curve and $\Z_4=\{1,i,i^2,i^3\}$ act on $E$ by multiplication. Let $\bL_0$ be an embedded circle of homology class $(1,1)$ and $\bL_1=i\bL_0$, $\bL_2=-\bL_0$, $\bL_3=-i\bL_0$ as Figure \ref{fig:244}. 
Then $\WT{\bL}=\bL_0\oplus\bL_1\oplus\bL_2\oplus\bL_3$ is weakly unobstructed with potential $W_{244}$. Every technical detail involved in the computation is just the same as $\Z_3$-case, so we only note that (Poincar\'e dual of) the oriented intersection $\WT{\bL} \cap C_h$ is $a-b-c+d$, and it is cohomologous to $p_3-p_0-p_1+p_2+{\rm (lower)}$ in $(CF(\WT{\bL},\WT{\bL}),m_1^{\WT{b}})$. For $C_v$, we obtain a cycle cohomologous to $p_0-p_1-p_2+p_3+{\rm (lower)}$. 
If we let $\widehat{\Z_4}=\{1,\chi,\chi^2,\chi^3\}$ such that $\chi(i)=i$, then 
\begin{equation}\label{eq:ks244} \ks_{orb}(PD[C_h])=\frac{-1-i}{2c_\bL}\xi_\chi+\frac{-1+i}{2c_\bL}\xi_{\chi^3}, \quad
\ks_{orb}(PD[C_v])=\frac{1-i}{2c_\bL}\xi_\chi+\frac{1+i}{2c_\bL}\xi_{\chi^3}.
\end{equation}

\subsubsection{$\Z_6$}
\begin{figure}
    \centering
    \includegraphics[width=0.85\linewidth]{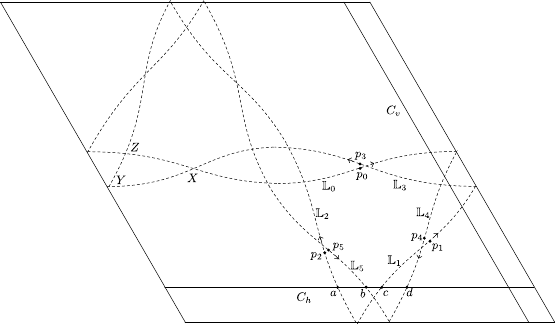}
    \caption{$\Z_6$-equivariant Lagrangian}
    \label{fig:236}
\end{figure}

Let $E=\C/(\Z+e^{2\pi i/3}\Z)$ and $\Z_6=\{1,\zeta=e^{\pi i/3},\zeta^2,\zeta^3,\zeta^4,\zeta^5\}$ act on $E$ by multiplication. Let $\bL_0$ be an embedded circle of class $(1,0)$ and $\bL_k=\zeta^k\bL_0$ for $k=1,\cdots,5$ as in Figure \ref{fig:236}. 
Let $\WT{\bL}:=\bigoplus_{k=0}^5 \zeta^k \bL_0$ be a weakly unobstructed Lagrangian with potential $W_{236}$. With the same computation as above, we obtain 
\[\CO^0(PD[C_h])\sim -p_2+p_5-p_1+p_4+{\rm (lower)}, \quad \CO^0(PD[C_v]) \sim p_1-p_3-p_4+p_0+{\rm (lower)},\]
hence for $\widehat{\Z_6}=\{1,\chi,\cdots,\chi^5\}$ with $\chi(\zeta)=\zeta$,
\begin{equation}\label{eq:ks236}
    \ks_{orb}(PD[C_h])=\frac{-\zeta-\zeta^2}{3c_\bL}\xi_\chi+\frac{\zeta+\zeta^2}{3c_\bL}\xi_{\chi^5},
    \quad
    \ks_{orb}(PD[C_v])=\frac{1+\zeta}{3c_\bL}\xi_\chi+\frac{1-\zeta^2}{3c_\bL}\xi_{\chi^5}.
\end{equation}


\section{Products in orbifold Jacobian algebras}\label{sec:orbiproduct}
Recall that our goal is to prove 
\[\ks_{orb}(c_\bL^2\cdot\pt_E) = -\frac{\det \; \Hess(W)}{|G|\cdot \mu}.\]
For this, we need to compute the product in the orbifold Jacobian algebra. Among various works, we follow the construction of \cite{LeeTwJac}. Other works on orbifold Jacobian algebras include \cite{BTW,ShkLGorb,HLL} etc.

Let $W\in \Bbbk[x_1,\cdots,x_n]$, and let $(W,\HG)$ be a Landau-Ginzburg orbifold. For $h\in \HG$, let 
\[ I_h:= \{i\mid h\cdot x_i \neq x_i\}\subset \{1,\cdots,n\},\quad  I^h:=I_h^c.\]
Let $S=\Bbbk[x_1,\cdots,x_n,x_1',\cdots,x_n']$. For $0\leq j<i\leq n$, define 
\[ \barW^h_{j,i}:=W(x_1',\cdots,x_{j}',x_{j+1},\cdots,x_i,hx_{i+1},\cdots,hx_n)\in S,\]
\[ \tildeW^h_{j,i}:=W(x_1^h,\cdots,x_{j}^h,x_{j+1},\cdots,x_{i},hx_{i+1},\cdots,hx_n)\in R.\]
(we define $x_i^h=x_i$ if $hx_i=x_i$ and $x_i^h=0$ otherwise.) We also define
\begin{equation}\label{eq:Wbar} \barW^h_{i,i}:=W(x_1',\cdots,x_{i}',hx_{i+1},\cdots,hx_n),\end{equation}
\begin{equation}\label{eq:Wtilde} \tildeW^h_{i,i}:=W(x_1^h,\cdots,x_i^h,hx_{i+1},\cdots,hx_n).\end{equation}
For $i,j\in \{1,\cdots,n\}$ with $j<i$, define
\begin{align}\label{eq:gjifji}
\begin{split}
g^h_{ji}&:=\begin{cases}
\frac{(\barW^h_{j,i}-\barW^h_{j-1,i})-(\barW^h_{j,i-1}-\barW^h_{j-1,i-1})}{(x_j'-x_j) (x_i-hx_i)} &{\rm if\;} i\in I_h, \\
0 & {\rm otherwise}, 
\end{cases} \\
 f^h_{ji}&:=
 \begin{cases}\frac{(\tildeW^h_{j,i}-\tildeW^h_{j-1,i})-(\tildeW^h_{j,i-1}-\tildeW^h_{j-1,i-1})}{(x_j-hx_j)(x_i-hx_i)}& {\rm if\;} i,j\in I_h,\\
 0 & {\rm otherwise},
 \end{cases}
 \end{split}
 \end{align}
 and
 \begin{align}\label{eq:gii}
 \begin{split}
      g^h_{ii}&:=
 \begin{cases}\frac{1}{x_i'-x_i}\cdot\Big(\frac{\barW^h_{i,i}-\barW^h_{i-1,i-1}}{x_i'-hx_i}-\frac{\barW^h_{i-1,i}-\barW^h_{i-1,i-1}}{x_i-hx_i}\Big) & {\rm if\;} i\in I_h, \\
0 & {\rm otherwise}.
 \end{cases}
 \end{split}\end{align}
Let $\theta_i$ and $\partial_i$ (for $i=1,\cdots,n$) be formal variables with $|\theta_i|=-1$, $|\partial_i|=1$ and
\begin{equation*}
 \theta_i \theta_j=-\theta_j \theta_i,\;\; \partial_{i}\partial_{j}=-\partial_{j}\partial_{i},\;\; \partial_{i}\theta_j=-\theta_j\partial_{i}+\delta_{ij}.
 \end{equation*}
For an ordered subset $I=\{i_1,\cdots,i_k\}\subset\{1,\cdots,n\}$, we introduce a notation
$ \theta_I:=\theta_{i_1}\cdots\theta_{i_k}.$ Now, define an $S$-linear map
\begin{align}
 \eta_h: S\langle\theta_1,\cdots,\theta_n,\partial_1,\cdots,\partial_n\rangle & \to S\langle\theta_1,\cdots,\theta_n,\partial_1,\cdots,\partial_n\rangle,\nonumber \\
\theta_I \partial_J& \mapsto  \sum (-1)^{|I|} g^h_{ji} \frac{\partial\theta_I}{\partial\theta_i}\partial_j\partial_J+\sum f^h_{ji}\frac{\partial^2 \theta_I}{\partial\theta_j\partial\theta_i}\partial_J. \label{etasign}
\end{align}
We also define a map
\begin{equation}\label{eq:exp} 
\exp(\eta_h):=1+\eta_h+\frac{\eta_h^2}{2!}+\cdots\;\;:S\langle\theta_1,\cdots,\theta_n,\partial_1,\cdots,\partial_n\rangle  \to S\langle\theta_1,\cdots,\theta_n,\partial_1,\cdots,\partial_n\rangle.
\end{equation}
For $f(x_1,\cdots,x_n,x_1',\cdots,x_n')\theta_I\partial_J \in S\langle \theta_1,\cdots,\theta_n,\partial_1,\cdots,\partial_n\rangle$, 
define 
\[ h_*(f\cdot\theta_I\partial_J):= f(x_1,\cdots,x_n,h^{-1}x_1',\cdots,h^{-1}x_n')\cdot h^{-1}(\theta_I\partial_J)\]
where for $h_i\in \Bbbk^*$ such that $h\cdot x_i=h_i x_i$, $h\cdot\theta_i = h_i^{-1}\theta_i$ and $h\cdot\partial_i=h_i\partial_i$.
Then for $h,h'\in \HG$, define
\[ \WT{\sigma}_{h,h'}:= \langle h'_*(\exp(\eta_{h})(\theta_{I_{h}}))\cdot \exp(\eta_{h'})(\theta_{I_{h'}}),\theta_{I_{hh'}}\rangle\]
which is $\theta_{I_{hh'}}$-coefficient of $h'_*(\exp(\eta_{h})(\theta_{I_{h}}))\cdot \exp(\eta_{h'})(\theta_{I_{h'}})$. Let 
\[ \pi_{h}: S \to S/(y_1-hx_1,\cdots,y_n-hx_n)\]
be the "$h$-twisted" quotient map for $h\in H$. Then
\[ \sigma_{h,h'}:= \pi_{hh'}(\WT{\sigma}_{h,h'}).\]
induces an element of $\Jac(W^{hh'})$, and is the structure constant of $\xi_h\bullet \xi_{h'}$, i.e.
\[ \xi_h\bullet \xi_{h'} = \sigma_{h,h'}\cdot\xi_{hh'}.\]

We present relevant multiplications inside three orbifold Jacobian algebras from elliptic curves.
\subsubsection{$(\C/(\Z+e^{2\pi i/3}\Z),\Z_3)$}
The mirror superpotential is given by
\[ W= \phi(x^3+y^3+z^3)-\psi xyz \in \Lambda[x,y,z]\]
and $\Z_3=\{1,\chi,\chi^2\}$ acts on $\Lambda[x,y,z]$ by
\[ \chi\cdot x= e^{2\pi i/3}x,\; \chi\cdot y= e^{2\pi i/3}y,\; \chi\cdot z= e^{2\pi i/3}z.\]
By $I_\chi=\{1,2,3\}=I_{\chi^2}$ (we let $x_1=x,\; x_2=y,\; x_3=z$) we deduce that $\xi_\chi$ and $\xi_{\chi^2}$ are $\Z_3$-invariant, hence are elements of $\Jac(W,\Z_3)$.

To calculate $\xi_{\chi}\bullet \xi_{\chi^2}$, let us first compute $\exp(\eta_\chi)(\theta_1\theta_2\theta_3)$ as follows. (Let $\rho=e^{2\pi i/3}$ for simplicity.)

\begin{align*}
    \barW^\chi_{0,0}&= W(x,y,z), &
    &\barW^\chi_{0,1}= \phi(x^3+y^3+z^3)-\rho^2 \psi xyz, \\
    \barW^\chi_{0,2}&=\phi(x^3+y^3+z^3)-\rho \psi xyz,&
    & \barW^\chi_{0,3}=W(x,y,z), \\
    \barW^\chi_{1,1}&= \phi(x'^3+y^3+z^3)-\rho^2 \psi x'yz, &
    &\barW^\chi_{1,2}= \phi(x'^3+y^3+z^3)-\rho \psi x'yz, \\
    \barW^\chi_{1,3}&=\phi(x'^3+y^3+z^3)- \psi x'yz,&
    & \barW^\chi_{2,2}=\phi(x'^3+y'^3+z^3)-\rho \psi x'y'z, \\
    \barW^\chi_{2,3}&= \phi(x'^3+y'^3+z^3)-  \psi x'y'z, &
    &\barW^\chi_{3,3}=W(x',y',z'),
\end{align*}
and
\begin{align*}
    \tildeW^\chi_{0,0}&= W(x,y,z), & &\tildeW^\chi_{0,1}= \phi(x^3+y^3+z^3)-\rho^2 \psi xyz, \\    
    \tildeW^\chi_{0,2}&=\phi(x^3+y^3+z^3)-\rho \psi xyz, & & \tildeW^\chi_{0,3}=W(x,y,z), \\
    \tildeW^\chi_{1,1}&=\tildeW^\chi_{1,2}=\tildeW^\chi_{1,3}=\phi(y^3+z^3), & & \tildeW^\chi_{2,2}=\tildeW^\chi_{2,3}=\phi z^3,\\
      \tildeW^\chi_{3,3}&=0. & &
\end{align*}
Plugging $\barW^\chi_{j,i}$ and $\tildeW^\chi_{j,i}$ into \eqref{eq:gjifji} and \eqref{eq:gii}, we have
\begin{align*}
    & g^\chi_{11}=\phi(x'-\rho^2 x),\quad g^\chi_{12}=-\rho\psi z,\quad g^\chi_{13}=-\psi y, \\ 
    & g^\chi_{22}=\phi(y'-\rho^2 y),\quad g^\chi_{23}=-\psi x',\quad g^\chi_{33}=\phi(z'-\rho^2 z), \\
    & f^\chi_{12}=\frac{\rho\psi z}{1-\rho}, \quad f^\chi_{13}=\frac{\psi y}{1-\rho},\quad f^\chi_{23}=0.
\end{align*}
By definitions \eqref{etasign} and \eqref{eq:exp},
\begin{align*}
    \exp(\eta_\chi)(\theta_1\theta_2\theta_3)=& \sum_{I \neq \emptyset}\sum_J f_{IJ}\theta_I\partial_J \\
    &+(g_{11}^\chi f_{23}^\chi  - g_{12}^\chi f_{13}^\chi +g_{13}^\chi f_{12}^\chi)\partial_1 + (-g_{22}^\chi f_{13}^\chi +g_{23}^\chi f_{12}^\chi)\partial_2 + g_{33}^\chi f_{12}^\chi \partial_3 \\
    &+ g_{11}^\chi g_{22}^\chi g_{33}^\chi \partial_3 \partial_2\partial_1 \\
    =& \sum_{I \neq \emptyset}\sum_J f_{IJ}\theta_I\partial_J \\
    & - \frac{\phi\psi(y'y-\rho^2 y^2)+\rho\psi^2 x'z}{1-\rho}\partial_2 + \frac{\rho\phi\psi(z'z-\rho^2 z^2)}{1-\rho}\partial_3\\
    &+\phi^3 (x'-\rho^2 x)(y'-\rho^2 y)(z'-\rho^2 z) \partial_3\partial_2\partial_1
\end{align*} 
and
\begin{align*}
   \pi_{\chi\chi^2}\big( \chi^2_* \exp(\eta_\chi)(\theta_1\theta_2\theta_3)\big)
   =& \chi^2_* \exp(\eta_\chi)(\theta_1\theta_2\theta_3)|_{x'=x,y'=y,z'=z} \\
   =& \chi^2_* \sum_{I \neq \emptyset}\sum_J f_{IJ}\theta_I\partial_J|_{x'=x,y'=y,z'=z} \\
   & - \frac{\phi\psi(\rho^2-1) y^2+\psi^2 xz}{1-\rho}\partial_2 + \frac{\phi\psi(1-\rho) z^2}{1-\rho}\partial_3\\
    &+3\phi^3 (\rho^2-\rho) xyz \cdot \partial_3\partial_2\partial_1.
\end{align*}
We do not have to compute $\sum_{I \neq \emptyset}\sum_J f_{IJ}\theta_I\partial_J$ explicitly, because it does not contribute to $\WT{\sigma}_{\chi,\chi^2}$ which is the coefficient of $\theta_{I_{\chi\chi^2}}=\theta_{I_1}=\theta_\emptyset$. 

The same computations as above gives $\exp(\eta_{\chi^2})(\theta_1\theta_2\theta_3)$ as follows.
\begin{align*}
    \exp(\eta_{\chi^2})(\theta_1\theta_2\theta_3) =& \sum_{I}\sum_{J\neq \emptyset}f_{IJ}\theta_I\partial_J + \theta_1\theta_2\theta_3 - f^{\chi^2}_{12}\theta_3+f^{\chi^2}_{13}\theta_2-f^{\chi^2}_{23}\theta_1 \\
    =& \sum_{I}\sum_{J\neq \emptyset}f_{IJ}\theta_I\partial_J+\theta_1\theta_2\theta_3-\frac{\rho^2 \psi z}{1-\rho^2}\theta_3+\frac{\psi y}{1-\rho^2}\theta_2.
\end{align*}
This time, we do not have to compute $\sum_{I}\sum_{J\neq \emptyset}f_{IJ}\theta_I\partial_J$ explicitly because it does not contribute to $\WT{\sigma}_{\chi,\chi^3}$. We conclude that
\begin{align}\label{eq:z3product}
\begin{split}
    {\sigma}_{\chi,\chi^3} =& \big(3\phi^3 (\rho^2-\rho)-\frac{\psi^3 }{3}\big)xyz
    -\frac{\phi\psi^2 (\rho^2-1)}{3}y^3 - \frac{\phi\psi^2(\rho^2-1)}{3}z^3.
\end{split}
\end{align}

\subsubsection{$(\C/(\Z+i\Z),\Z_4)$}
The mirror superpotential is given by
\[ W=-qxyz+q^6 x^2 + a(y^4+z^4)+by^2 z^2 \in \Lambda[x,y,z]\]
and $\Z_4=\{1,\chi,\chi^2,\chi^3\}$ acts on $\Lambda[x,y,z]$ by
\[ \chi \cdot x= -x, \; \chi\cdot y= iy, \; \chi\cdot z=iz.\]
By $I_\chi=\{1,2,3\}=I_{\chi^3}$ (we let $x_1=x,\; x_2=y,\; x_3=z$) we deduce that $\xi_\chi$ and $\xi_{\chi^3}$ are $\Z_4$-invariant, hence are elements of $\Jac(W,\Z_4)$.

As above we compute $\exp(\eta_\chi)(\theta_1\theta_2\theta_3)$ from
\begin{align*}
    \barW^\chi_{0,0}&= W(x,y,z), &
    &\barW^\chi_{0,1}= q^6 x^2+qxyz+ay^4+az^4+by^2 z^2, \\
    \barW^\chi_{0,2}&=q^6 x^2-qixyz+ay^4+az^4-by^2 z^2,&
    & \barW^\chi_{0,3}=W(x,y,z), \\
    \barW^\chi_{1,1}&= q^6 x'^2+qx'yz+ay^4+az^4+by^2 z^2, &
    &\barW^\chi_{1,2}= q^6 x'^2-qixyz+ay^4+az^4-by^2 z^2, \\
    \barW^\chi_{1,3}&=q^6 x'^2-qx'yz+ay^4+az^4+by^2 z^2,&
    & \barW^\chi_{2,2}=q^6 x'^2-qix'y'z+ay'^4+az^4-by'^2 z^2, \\
    \barW^\chi_{2,3}&= q^6 x'^2-qx'y'z+ay'^4+az^4+by'^2 z^2, &
    &\barW^\chi_{3,3}=W(x',y',z'),
\end{align*}
and
\begin{align*}
    \tildeW^\chi_{0,0}&= W(x,y,z), & &\tildeW^\chi_{0,1}= q^6 x^2+qxyz+ay^4+az^4+by^2 z^2, \\    
    \tildeW^\chi_{0,2}&=q^6 x^2-qixyz+ay^4+az^4-by^2 z^2, & & \tildeW^\chi_{0,3}=W(x,y,z), \\
    \tildeW^\chi_{1,1}&=ay^4+az^4+by^2 z^2, & & \tildeW^\chi_{1,2}=ay^4+az^4-by^2 z^2, \\
    \tildeW^\chi_{1,3}&=ay^4+az^4+by^2 z^2, & &\tildeW^\chi_{2,2}=\tildeW^\chi_{2,3}=az^4,\\
      \tildeW^\chi_{3,3}&=0. & &
\end{align*}
Plugging $\barW^\chi_{j,i}$ and $\tildeW^\chi_{j,i}$ into \eqref{eq:gjifji} and \eqref{eq:gii}, we have
\begin{align*}
    & g^\chi_{11}=q^6,\quad g^\chi_{12}=-qiz,\quad g^\chi_{13}=-qy, \\ 
    & g^\chi_{22}=a(y'+y)(y'+iy)-bz^2,\quad g^\chi_{23}=-qx'+b(y'+y)(z+iz),\quad g^\chi_{33}=a(z'+z)(z'+iz)+by'^2, \\
    & f^\chi_{12}=-\frac{qiz}{2}, \quad f^\chi_{13}=\frac{qy}{2},\quad f^\chi_{23}=-ibyz.
\end{align*}
By definitions \eqref{etasign} and \eqref{eq:exp},
\begin{align*}
    \exp(\eta_\chi)(\theta_1\theta_2\theta_3)=& \sum_{I \neq \emptyset}\sum_J f_{IJ}\theta_I\partial_J \\
    &+(g_{11}^\chi f_{23}^\chi  - g_{12}^\chi f_{13}^\chi +g_{13}^\chi f_{12}^\chi)\partial_1 + (-g_{22}^\chi f_{13}^\chi +g_{23}^\chi f_{12}^\chi)\partial_2 + g_{33}^\chi f_{12}^\chi \partial_3 \\
    &+ g_{11}^\chi g_{22}^\chi g_{33}^\chi \partial_3 \partial_2\partial_1 \\
    =& \sum_{I \neq \emptyset}\sum_J f_{IJ}\theta_I\partial_J \\
    & -q^6 biyz \partial_1 + \frac{1}{2}\big( qbyz^2-qa(y'+y)(y'+iy)y+(-qx'+b(y'+y)(z+iz))qiz\big)\partial_2 \\
    &+ (a(z'+iz)(z'+z)+by'^2)\cdot \frac{qiz}{2}\partial_3 \\
    &+ q^6(a(y'+y)(y'+iy)-bz^2)(a(z'+z)(z'+iz)+by'^2)\partial_3 \partial_2 \partial_1,
\end{align*} 
and
\begin{align*}
   \pi_{\chi\chi^3}\big( \chi^3_* \exp(\eta_\chi)(\theta_1\theta_2\theta_3)\big)
   =& \chi^3_* \exp(\eta_\chi)(\theta_1\theta_2\theta_3)|_{x'=x,y'=y,z'=z} \\
   =& \chi^3_* \sum_{I \neq \emptyset}\sum_J f_{IJ}\theta_I\partial_J|_{x'=x,y'=y,z'=z} \\
   &+q^6 biyz \partial_1 + \big(-qbyz^2-2i(i+1)qay^3+q^2ixz-2(i+1)qbyz^2\big)\frac{i}{2}\partial_2\\
   &+\big(-qa(i+1)z^3-\frac{qiby^2 z}{2}\big)i\partial_3\\
   &+q^6((2i-2)ay^2-bz^2)((2i-2)az^2-by^2)\partial_3\partial_2\partial_1
\end{align*}
We do not have to compute $\sum_{I \neq \emptyset}\sum_J f_{IJ}\theta_I\partial_J$ explicitly, because it does not contribute to $\WT{\sigma}_{\chi,\chi^3}$ which is the coefficient of $\theta_{I_{\chi\chi^3}}=\theta_{I_1}=\theta_\emptyset$. 

The same computations as above gives $\exp(\eta_{\chi^3})(\theta_1\theta_2\theta_3)$ as follows.
\begin{align*}
    \exp(\eta_{\chi^3})(\theta_1\theta_2\theta_3) =& \sum_{I}\sum_{J\neq \emptyset}f_{IJ}\theta_I\partial_J + \theta_1\theta_2\theta_3 - f^{\chi^3}_{12}\theta_3+f^{\chi^3}_{13}\theta_2-f^{\chi^3}_{23}\theta_1 \\
    =& \sum_{I}\sum_{J\neq \emptyset}f_{IJ}\theta_I\partial_J+\theta_1\theta_2\theta_3+\frac{qiz}{2}\theta_3+\frac{qy}{2}\theta_2-biyz\theta_1.
\end{align*}
This time, we do not have to compute $\sum_{I}\sum_{J\neq \emptyset}f_{IJ}\theta_I\partial_J$ explicitly because it does not contribute to $\WT{\sigma}_{\chi,\chi^3}$. We conclude that
\begin{align}\label{eq:z4product}
\begin{split}
    {\sigma}_{\chi,\chi^3} =&-\frac{q^3 xyz}{4}+\big(\frac{q^2a(i+1)}{2}-4q^6i(i+1)ab\big)y^4+ \frac{q^2 a(i+1)z^4}{2}  \\
    &+ \big(2q^6 b^2-8q^2 a^2 i -\frac{q^2 b(i-1)}{2}\big) y^2 z^2.
\end{split}
\end{align}

\subsubsection{$(\C/(\Z+e^{2\pi i/3}\Z),\Z_6)$}
The mirror superpotential is given by
\[ W=-qxyz+q^6 x^2 + a_1 y^3 + a_2 z^6 + a_3 y^2 z^2 + a_4 yz^4\in \Lambda[x,y,z]\]
and $\Z_6=\{1,\chi,\chi^2,\chi^3,\chi^4,\chi^5\}$ acts on $\Lambda[x,y,z]$ by
\[ \chi\cdot x= -x,\; \chi\cdot y= e^{2\pi i/3}y, \; \chi\cdot z = e^{\pi i/3}z.\]
By $I_\chi=\{1,2,3\}=I_{\chi^5}$ we deduce that $\xi_\chi$ and $\xi_{\chi^5}$ are $\Z_6$-invariant, so they are elements of $\Jac(W,\Z_6)$. Let $\zeta=e^{\pi i/3}$ for simplicity. To compute $\xi_\chi \bullet \xi_{\chi^5}$, we need
\begin{align*}
    \barW^\chi_{0,0}&= W(x,y,z), \\
    \barW^\chi_{0,1}&= q^6 x^2+qxyz+a_1 y^3 +a_2 z^6 +a_3 y^2 z^2+ a_4 yz^4, \\
    \barW^\chi_{0,2}&=q^6 x^2-\zeta qxyz+a_1 y^3+a_2 z^6 + \zeta^2 a_3 y^2 z^2+ a_4 yz^4, \\ \barW^\chi_{0,3}&=W(x,y,z), \\
    \barW^\chi_{1,1}&= q^6 x'^2+qx'yz+a_1 y^3 +a_2 z^6 +a_3 y^2 z^2+ a_4 yz^4,  \\
    \barW^\chi_{1,2}&= q^6 x'^2-\zeta qx'yz+a_1 y^3 +a_2 z^6 +\zeta^2 a_3 y^2 z^2+ \zeta^4 a_4 yz^4, \\
    \barW^\chi_{1,3}&=q^6 x'^2-qx'yz+a_1 y^3 +a_2 z^6 +a_3 y^2 z^2+ a_4 yz^4, \\
    \barW^\chi_{2,2}&=q^6 x'^2-\zeta qx'y'z+a_1 y'^3 +a_2 z^6 +\zeta^2 a_3 y'^2 z^2+ \zeta^4 a_4 y'z^4, \\
    \barW^\chi_{2,3}&= q^6 x'^2+qx'y'z+a_1 y'^3 +a_2 z^6 +a_3 y'^2 z^2+ a_4 y'z^4, \\
    \barW^\chi_{3,3}&=W(x',y',z'),
\end{align*}
and
\begin{align*}
    \tildeW^\chi_{0,0}&= W(x,y,z), \\
    \tildeW^\chi_{0,1}&= q^6 x^2+qxyz+a_1 y^3 +a_2 z^6 +a_3 y^2 z^2+ a_4 yz^4, \\    
    \tildeW^\chi_{0,2}&=q^6 x^2-\zeta qxyz+a_1 y^3 +a_2 z^6 + \zeta^2 a_3 y^2 z^2+ \zeta^4 a_4 yz^4, \\
    \tildeW^\chi_{0,3}&=W(x,y,z), \\
    \tildeW^\chi_{1,1}&=a_1 y^3 +a_2 z^6 +a_3 y^2 z^2+ a_4 yz^4, \\
    \tildeW^\chi_{1,2}&=a_1 y^3 +a_2 z^6 +\zeta^2 a_3 y^2 z^2+\zeta^4 a_4 yz^4, \\
    \tildeW^\chi_{1,3}&= a_1 y^3 + a_2 z^6+ a_3 y^2 z^2 + a_4 yz^4, \\
    \tildeW^\chi_{2,2}&=\tildeW^\chi_{2,3}=az^4, \\
    \tildeW^\chi_{3,3}&=0. 
\end{align*}
As above, we calculate the following:
\begin{align*}
    & g^\chi_{11}=q^6,\quad g^\chi_{12}=-\zeta qz,\quad g^\chi_{13}=-qy, \\ 
    & g^\chi_{22}=a_1 (y'+\zeta y)+\zeta^2 a_3 z^2,\quad g^\chi_{23}=-qx'+a_3(y'+y)(z+\zeta z)+a_4 \sqrt{3}iz^3,\\
    & g^\chi_{33}=a_2(z'-\zeta^2 z)(z'-\zeta^3 z)(z'-\zeta^4 z)(z'-\zeta^5 z)+a_3 y'^2
    +a_4 y'(z'^2+z'z+z^2+(z'+z)\zeta z+\zeta^2 z^2), \\
    & f^\chi_{12}=-\frac{\zeta qz}{2}, \quad f^\chi_{13}=\frac{qy}{2},\quad f^\chi_{23}=\frac{-a_3 yz-\zeta a_4 z^3}{1-\zeta}
\end{align*}
and 
\begin{align*}
    \pi_{\chi\chi^5}\big( \chi^5_* \exp(\eta_\chi)(\theta_1\theta_2\theta_3)\big) =& \chi^5_*\sum_{I\neq \emptyset}\sum_J f_{IJ}\theta_I\partial_J |_{x'=x,y'=y,z'=z} \\
    &+q^6\cdot \frac{a_3 yz+\zeta a_4 z^3}{1-\zeta}\partial_1 \\
    &+\big( (-a_1\sqrt{3}iy-\zeta^2 a_3 z^2)\frac{qy}{2}+(qx+a_3 \sqrt{3}iyz+a_4 \sqrt{3}iz^3)\frac{\zeta qz}{2}\big) \zeta^2 \partial_2 \\
    &+ \big( -6a_2 z^4+a_3 \zeta^4 y^2 + a_4(3\zeta^4-1)yz^2\big)\frac{\zeta^2 qz}{2}\partial_3 \\
    &+ q^6 ( -6a_2 z^4+a_3 \zeta^4 y^2 + a_4(3\zeta^4-1)yz^2)(a_1 \sqrt{3}iy+\zeta^2 a_3 z^2)\partial_3\partial_2\partial_1.
\end{align*}
Also, we have
\[ \exp(\eta_{\chi^5})(\theta_1\theta_2\theta_3)= \sum_I\sum_{J\neq \emptyset} f_{IJ}\theta_I\partial_J 
+\theta_1\theta_2\theta_3 -\frac{\zeta^{-1}qz}{2}\theta_3+\frac{qy}{2}\theta_2+\frac{a_3 yz+\zeta^{-1}a_4 z^3}{1-\zeta^{-1}}\theta_1.\]
In conclusion, 
\begin{align*}
    \sigma_{\chi,\chi^5} 
    =& \big( -\frac{a_1 \sqrt{3}i q^2 \zeta^2}{4}+q^6 a_1 a_3 \zeta^4\sqrt{3}i\big)y^3 \\
    &+ \big( {q^6 a_4^2}+ \frac{3q^2 a_2\zeta}{2}-6q^6\zeta^2 a_2 a_3 \big)z^6 \\
    &+ \big( \frac{(9-\sqrt{3}i)q^6 a_3 a_4}{2} +\frac{(-1+\sqrt{3}i)q^2 a_4}{4}-6q^6 a_1 a_2 \sqrt{3}i\big)yz^4 \\
    &+ \big( 2q^6 a_3^2 + q^6 a_1 a_4(3\zeta^4-1)\sqrt{3}i\big) y^2 z^2 - \frac{q^3 xyz}{4}.
\end{align*}



\section{Main results}\label{sec:mainthm}
A {\em Frobenius algebra} is a unital associative algebra over $\Bbbk$ together with a nondegenerate bilinear form $\langle\cdot,\cdot\rangle$ such that $\langle xy,z\rangle=\langle x,yz\rangle$. Observe that given a nondegenerate pairing $\langle\cdot,\cdot\rangle$ on a Frobenius algebra, the rescaling 
$ \langle c\cdot, c\cdot\rangle$ for $c\neq 0$ can be still used to define a new Frobenius algebra with the same ring structure. Given two Frobenius algebras $(A,\langle\;,\;\rangle_A)$ and $(B,\langle\;,\;\rangle_B)$, $f:A \to B$ is an isomorphism of Frobenius algebras iff it is a ring isomorphism and an isometry with respect to nondegenerate pairings. Given unities $1_A\in A$ and $1_B\in B$, let ${\rm{tr}}_A:A \to \Bbbk$ and ${\rm{tr}}_B:B\to \Bbbk$ be maps (a.k.a. {\em trace}s) given by $\langle\cdot,1_A\rangle_A$ and $\langle\cdot,1_B\rangle_B$ respectively. Then $f:A \to B$ is an isometry if and only if it preserves traces.

For a compact symplectic manifold, the quantum cup product together with the Poincar\'e pairing $\langle \cdot,\cdot\rangle_{PD}$ gives rise to a Frobenius algebra structure on $QH^*(X)$. 
For an isolated singularity $W\in \Bbbk[x_1,\cdots,x_n]$, we equip $\Jac(W)$ with the residue pairing
\[ \langle f,g\rangle_{res}=(-1)^{\frac{n(n-1)}{2}}{\rm{Res}}\begin{bmatrix}
        fg\cdot dx_1\wedge\cdots\wedge dx_n \\ \partial_{x_1}W,\cdots,\partial_{x_n}W
    \end{bmatrix}.\]
As pointed out in \cite{CLS}, the Kodaira-Spencer map does not intertwine pairings $\langle\cdot,\cdot\rangle_{PD}$ and $\langle\cdot,\cdot\rangle_{res}$ in general. Examining the Cardy condition, we are led to consider a modification of the residue pairing by the constant given by ratio of Floer volume form and the usual volume form on the reference Lagrangian $\bL$. In Section \ref{sec:KSelliptic} we denoted such ratios by $c_\bL$. 
\begin{theorem}
Let $X$ be an elliptic orbisphere and $W\in \Lambda[x,y,z]$ be the LG mirror to $X$. Let $A=(QH^*(X),\langle\cdot,\cdot\rangle_{PD})$ and $B=(\Jac(W),\langle c_\bL\cdot,c_\bL\cdot\rangle_{res})$ be Frobenius algebras where $c_\bL$ is the rescaling constant of Seidel Lagrangian $\bL\subset X$. Then
\begin{equation}\label{eq:mainthm}
c_\bL^2\cdot \ks(\pt_X) = -\frac{\det \; {\rm Hess}(W)}{\mu},
\end{equation}
where $\mu=\dim \Jac(W)$.
As a result, $\ks: A \to B$ is an isomorphism of Frobenius algebras, i.e. it preserves traces as follows:
\begin{equation*} \langle \pt_X,1_X \rangle_{PD} = \langle c_\bL\cdot \ks(\pt_X),c_\bL \cdot 1 \rangle_{res}.
\end{equation*}
\end{theorem}
Provided \eqref{eq:mainthm}, the second statement follows immediately from a classical result on the residue over an isolated singularity (the author thanks Cheol-Hyun Cho for pointing out this fact).
\begin{theorem}[Section 17, \cite{HaRes}]
\[{\rm{Res}}\begin{bmatrix}
        {\det\;\rm{Hess}}(W)\cdot dx_1\wedge\cdots\wedge dx_n \\ \partial_{x_1}W,\cdots,\partial_{x_n}W
    \end{bmatrix}=\mu.\]
\end{theorem}
We can easily compare $\ks:QH^*(X) \to \Jac(W)$ and $\ks_{orb}: H^*(E) \to \Jac(W,\HG)$ as follows. Basically, holomorphic discs employed in the computation of $\ks_{orb}(\pt_{E})$ are same as those used for $\ks(\pt_X)$, but the discs which are related by $G$-action on $E$ are identified as the same disc on $X$. Thus, $|G|$ different discs which are used for $\ks_{orb}(\pt_E)$ project to a single disc for $\ks(\pt_X)$. We conclude that
\[ \ks_{orb}(\pt_E)= \frac{1}{|G|}\cdot \ks(\pt_X).\]
The rest of the paper is thus devoted to prove
\[ c_\bL^2\cdot \ks_{orb}(\pt_E)=-\frac{\det \; {\rm Hess}(W)}{|G|\cdot \mu}\]
where $G=\Z_3,\Z_4$ or $\Z_6$ acts on $E$, with $W$ given accordingly with respect to $G$-equivariant Seidel Lagrangian.

\subsection{$\Z_3$-case}
Recall from \eqref{eq:ks333} that (with $\rho=e^{2\pi i/3}$)
\[ c_\bL\cdot \ks_{orb}(PD[C_h])=-\rho \xi_\chi-\rho^2 \xi_{\chi^2},\quad 
c_\bL\cdot \ks_{orb}(PD[C_v])= -\rho^2 \xi_{\chi}-\rho \xi_{\chi^2},\]
hence by \eqref{eq:z3product}
\begin{align*} c_\bL^2\cdot \ks_{orb}(\pt_E)&=c_\bL^2 \ks_{orb}(PD[C_h]) \bullet \ks_{orb}(PD[C_v]) \\
&= (\rho^2-\rho)\xi_\chi \bullet \xi_{\chi^2} \\
&= (\rho^2-\rho)\cdot \Big( \big(3\phi^3 (\rho^2-\rho)-\frac{\psi^3 }{3}\big)xyz
    -\frac{\phi\psi^2 (\rho^2-1)}{3}y^3 - \frac{\phi\psi^2(\rho^2-1)}{3}z^3\Big).
\end{align*}
Observe that we used supercommutativity of the product for the second identity. The Jacobian ideal $\partial W$ is given by
\[ \partial W = (3\phi x^2-\psi yz, 3\phi y^2-\psi xz, 3\phi z^2-\psi xy),\]
so $y^3=z^3=\frac{\psi xyz}{3\phi}$ modulo $\partial W$, and
\[ c_\bL^2\cdot \ks_{orb}(\pt_E) = \big(-9\phi^3 + \frac{\psi^3}{3}\big) xyz.\]
On the other hand,
\[ {\rm Hess}(W)=\begin{pmatrix}
    6\phi x & -\psi x & -\psi y \\ -\psi z & 6\phi y & -\psi x \\ -\psi y & -\psi x & 6\phi z
\end{pmatrix},\]
hence $\det\; {\rm Hess}(W)= (216\phi^3-8\psi ^3) xyz$ in $\Jac(W)$. By $\mu=8$ and $|G|=3$, we easily verify \eqref{eq:mainthm}.

\subsection{$\Z_4$-case}
From \eqref{eq:ks244} we have
\[ c_\bL^2\cdot \ks_{orb}(\pt_E)=-i\xi_\chi\bullet \xi_{\chi^3}=-i\sigma_{\chi,\chi^3}\]
where $\sigma_{\chi,\chi^3}$ as \eqref{eq:z4product}. Modulo $\partial W= (2q^6 x-qyz, 4ay^3-qxz+2byz^2,4az^3-qxy+2by^2 z)$, we have relations
\[ xyz=2q^5 x^2,\quad y^2 z^2= 4q^{10}x^2,\quad 2ay^4=2az^4=(q^6-4bq^{10})x^2,\]
hence in $\Jac(W)$,
\begin{align}
-i \sigma_{\chi,\chi^3}&= -i \big(\frac{iq^8}{2}-4iq^{12}b+8iq^{16}b^2-32iq^{16}a^2\big)x^2\nonumber\\
&= \big(\frac{q^8}{2}-4q^{12}b+8q^{16}b^2-32q^{16}a^2\big)x^2.\nonumber
\end{align}
By
\[ {\rm Hess}(W)= \begin{pmatrix}
    2q^6 & -qz & -qy \\ -qz & 12ay^2+2bz^2 & -qx+4byz \\ -qy & -qx+4byz & 12az^2+2by^2
\end{pmatrix}\]
we have
\begin{align}
&\det\; {\rm Hess}(W)\nonumber\\
=& (48q^6 ab-12q^2 a)(y^4+z^4)-2q^8 x^2 + (288q^6 a^2-24q^6 b^2 + 4q^2 b)y^2 z^2 + (16q^7 b-2q^3)xyz \nonumber\\
=& (-18q^8+144q^{12}b-288q^{16}+1152q^{16}a^2)x^2. \nonumber
\end{align}
This time, $|G|=4$ and $\mu=9$, and it is straightforward that \eqref{eq:mainthm} holds.

\subsection{$\Z_6$-case}
From \eqref{eq:ks236} we have 
\[ c_\bL^2 \cdot \ks_{orb}(\pt_E)=-\frac{\sqrt{3}i}{3}\xi_{\chi}\bullet \xi_{\chi^5} 
=-\frac{\sqrt{3}i}{3} \sigma_{\chi,\chi^5}.\]
We also compute
\begin{align*}
    \det \; {\rm Hess}(W) =& (120q^6 a_2 a_3 - 32 q^6 a_4^2 -30q^2 a_2)z^6 + (360q^6 a_1 a_2-16q^6 a_3 a_4-4q^2 a_4)yz^4 \\
    &+(144q^6 a_1 a_4 -24q^6 a_3^2 + 4q^2 a_3)y^2 z^2 + 16q^7 a_4 xz^3 + (24q^6 a_1 a_3 - 6q^2 a_1)y^3 \\
    &+ (16q^7 a_3 - 2q^3 )xyz - 2q^8 x^2.
\end{align*}
The relations modulo $\partial W$ needed for our purpose are as follows:
\[ xyz= 2q^5 x^2, \quad a_1 y^3 = \frac{-8q^{14}a_3 a_4^2 - 16q^{10}a_2 a_3 + 32q^{14}a_2 a_3^2 -24q^{14}a_1a_2a_4 + 2q^{10}a_4^2 + 2q^6 a_2}{4q^4 a_4^2 + 3a_2-12q^4 a_2 a_3}x^2,\]
\[ z^6 = \frac{q^6(-8q^4 a_3 - 48q^8 a_1 a_4 + 16q^8 a_3^2 +1)}{4q^4 a_4^2 + 3a_2-12q^4 a_2 a_3}x^2,\quad
yz^4 = \frac{2q^{10}(36q^4 a_1 a_2 -4q^4 a_3 a_4 +a_4) }{4q^4 a_4^2 + 3a_2-12q^4 a_2 a_3}x^2.\]
The result \eqref{eq:mainthm} follows from $|G|=6$ and $\mu=10$. We omit the tedious computation which involves substitutions of above relations to $\det\; {\rm Hess}(W)$ and $\sigma_{\chi,\chi^5}$.

\subsection{Nontrivial identitites of formal power series}
Recall from \cite{ACHL} that \[\ks(\pt_X)=\frac{1}{8}q\cdot \frac{\partial W}{\partial q}\]
where $X$ is an orbisphere $\PP^1_{a,b,c}$ with area $1$ (If the area is $A$, then we take $\frac{1}{8A}$ instead of $\frac{1}{8}$). Hence we also have
\[ \ks_{orb}(\pt_E)=\frac{1}{8|G|}q\cdot \frac{\partial W}{\partial q}.\]
We computed $\ks_{orb}(\pt_E)$ in another way, namely by orbifold Jacobian algebra structure. Let $W_{333}$, $W_{244}$ and $W_{236}$ be mirror superpotentials from Seidel Lagrangians $\bL_{333}\subset \PP^1_{3,3,3}$, $\bL_{244}\subset\PP^1_{2,4,4}$ and $\bL_{236}\subset\PP^1_{2,3,6}$ respectively. 

Let us consider $\PP^1_{3,3,3}$. Recall that
\[ W_{333}= \phi(x^3+y^3+z^3)-\psi xyz\]
with its coefficients are 
\[ \phi= \sum_{k\in \Z} (-1)^{k+1} (k+\frac{1}{2})q^{(6k+3)^2},\quad \psi= \sum_{k\in \Z} (-1)^{k+1}(6k+1)q^{(6k+1)^2}.\]
We also have 
\[c_{\bL_{333}}= \sum_{k\in \Z} (-1)^k q^{(6k+1)^2}.\]
By the following
\begin{align*} \frac{c_{\bL_{333}}^2}{24}q\cdot \big(\frac{\partial \phi}{\partial q}(x^3+y^3+z^3)-\frac{\partial \psi}{\partial q}xyz\big) 
= \frac{c_{\bL_{333}}^2}{24}q\cdot \big( \frac{\psi}{\phi} \frac{\partial\phi}{\partial q}-\frac{\partial \psi}{\partial q}\big) xyz 
=\big(-9\phi^3 + \frac{\psi^3}{3}\big) xyz,
\end{align*}
we derive an identity of two arithmetics of formal power series
\[ \frac{c_{\bL_{333}}^2}{24}q\cdot \big( \frac{\psi}{\phi} \frac{\partial\phi}{\partial q}-\frac{\partial \psi}{\partial q}\big)= \big(-9\phi^3 + \frac{\psi^3}{3}\big).\]
Note that the identity was first proved in \cite{CLS} using theory of modular forms, and now we have a new proof without appealing to number theory.

The same argument can be applied to other two cases. For $W_{244}$, we have
\begin{align*}
    \frac{q^3}{32}\big( 6q^5 x^2 - xyz+ \frac{\partial a}{\partial q}(y^4+z^4)+\frac{\partial b}{\partial q}y^2 z^2 \big) &= \frac{q^3}{32}\big( 4q^5 +\frac{1}{a}\frac{\partial a}{\partial q}\cdot (q^6-4bq^{10}) + 4q^{10} \frac{\partial b}{\partial q}\big)x^2\\
    &= \big(\frac{q^8}{2}-4q^{12}b+8q^{16}b^2-32q^{16}a^2\big)x^2
\end{align*}
which implies
\[ \frac{q^3}{32}\big( 4q^5 +\frac{1}{a}\frac{\partial a}{\partial q}\cdot (q^6-4bq^{10}) + 4q^{10} \frac{\partial b}{\partial q}\big) = \frac{q^8}{2}-4q^{12}b+8q^{16}b^2-32q^{16}a^2\]
where $a$ and $b$ are formal power series given by
\[ a=\sum_{r\geq 0} (2r+1)q^{16(2r+1)^2-4}+\sum_{s>r\geq 0} (2r+2s+2)q^{16(2r+1)(2s+1)-4},\]
\[ b=\sum_{r\geq 1, s\geq 1}\big( -(4r+4s-2)q^{16(2r-1)2s-4}+(2r+2s)q^{64rs-4}\big).\]
We note that $c_{\bL_{2,4,4}}=\pm q$, so its square is $q^2$. For $W_{236}$ we can obtain a similar type of identity which is omitted.

\bibliographystyle{amsalpha}
\bibliography{geometry}

\end{document}